
\documentstyle{amsppt}
\magnification=\magstep 1
\TagsOnRight
\NoBlackBoxes
\leftheadtext{A\.G\. Leiderman, S\.A\. Morris and V\.G\. Pestov}
\rightheadtext{Free abelian topological group}

\def\R {{\Bbb R}}

\def\N{{\Bbb N}}

\topmatter
\title
The free abelian topological group
and the free locally convex space on the unit interval \thinspace\dag
\endtitle
\author
 Arkady G\. Leiderman, $^1$ Sidney A\. Morris \\
{\rm and}
Vladimir G\. Pestov $^2$
\endauthor
\affil
Ben-Gurion University of the Negev, P.O. Box 653,
Beer-Sheva, {\smc Israel} \\
University of Wollongong, Wollongong, N.S.W. 2522, {\smc Australia} \\
Victoria University of Wellington, P.O. Box 600, Wellington, {\smc New
Zealand}
\endaffil
\abstract{
We give a complete description of the topological spaces
$X$ such that the free abelian topological group
$A(X)$ embeds into the free abelian topological group $A(I)$
of the closed unit interval.
In particular, the free abelian topological group $A(X)$ of
any finite-dimensional compact metrizable space $X$ embeds into
$A(I)$.
The situation turns out to be somewhat
different for free locally convex spaces.
Some results for the spaces of continuous functions with the pointwise
topology are also obtained. Proofs are based on the classical
Kolmogorov's Superposition Theorem.}
\endabstract
\endsubjclass
\endkeywords
\endtopmatter
\document
\footnote""{\dag\thinspace Dedicated to the memory of Eli Katz.}
\footnote""{$^1$\thinspace Research Supported by the Israel Ministry of
Science.}
\footnote""{$^2$\thinspace Partially supported by
the Internal Grants Committee of the
Victoria University of Wellington.}
\smallpagebreak
\heading
\S 1. Introduction
\endheading
\smallpagebreak
The following natural question
arises as a part of the search for
a topologized version of the Nielsen-Schreier subgroup theorem. Let $X$
and $Y$ be completely regular
topological spaces; in which cases
the free (free abelian) topological group over $X$ can be embedded
as a
topological subgroup into the free (free abelian) topological group
over $Y$?
This problem has been treated for a long time
\cite{4, 10, 12-17, 21-23, 25, 28},
ever since it became clear that in general a topological subgroup of a
free (free abelian) topological group need not be topologically free
\cite{8, 4, 10}.
Recently a complete answer was obtained in the case where $X$ is a
subspace of
$Y$ and the embedding of free
topological groups extends the embedding of
spaces \cite{35}.
However, we are interested in the existence of an embedding which is
not necessarily a ``canonical'' one.
Among the most notable achievements, there are
certain sufficient conditions for a
subgroup of a free topological group to be
topologically free \cite{4, 22} and the following results.
\proclaim{Theorem 1.1 {\rm \cite{14}}}
If $X$ is a closed topological subspace of the free topological group
$F(I)$ then the free topological group $F(X)$ is a closed topological
subgroup of $F(I)$, where $I$ is the closed unit interval.
\qed\endproclaim
\proclaim{Corollary 1.2 {\rm \cite{22}}}
If $X$ is a finite-dimensional metrizable compact space then
$F(X)$ is a closed topological
subgroup of $F(I)$.
\qed\endproclaim
The abelian case proved to be more difficult, and
the following is the strongest result
known to date.
\proclaim{Theorem 1.3 \cite{12}}
If $X$ is a countable CW-complex of dimension $n$, then
the free abelian topological group on $X$ is a closed subgroup of the
free abelian topological group on the closed ball $B^n$.
\qed\endproclaim
\proclaim{Corollary 1.4 \cite{13}}
$A(\R)$ embeds into $A(I)$ as a closed topological subgroup.
\qed\endproclaim
It is known \cite{29}
that the covering dimension of any two free topological bases in a free
(abelian) topological group is the same;
this result is similar to the well-known
property of free bases of a discrete free (abelian) group having the
same cardinality, called the rank of the group.
Since the rank of a
subgroup of a free abelian group
cannot exceed the
rank of the group itself, it was conjectured \cite{15, 20} that the
dimension of a topological basis
of a topologically free
subgroup of a free
abelian topological group  $A(X)$ cannot exceed $dim~X$.
It remained even unclear
whether the group $A(I^2)$ embeds into $A(I)$
\cite{15}.

In this paper we prove that if $X$ is
a completely regular space then the free
abelian topological group $A(X)$
embeds into $A(I)$ as a topological subgroup if and only if
$X$ is a submetrizable $k_\omega$-space such that
every compact subspace of $X$ is
finite-dimensional. Another characterization: $X$
is homeomorphic to a closed topological subspace of the group $A(I)$
itself. In particular, if
$X$ is a compact metrizable space of finite dimension, then
$A(X)$ embeds into $A(I)$. Thus, the
analogy with the non-abelian case is complete.
We also study the
problem of embedding the free locally convex space $L(X)$ into
the free locally convex space $L(I)$ and characterize those
$k_\omega$-spaces $X$ admitting such an embedding.
Paradoxically, such spaces $X$ are just all compact metrizable
finite-dimensional spaces.
In particular, the free LCS $L(\R)$ does not embed
into $L(I)$.
Our results provide answers to a number of open problems from
\cite{20, 15, 27}.

\smallpagebreak
\heading
\S 2. Preliminaries
\endheading
\smallpagebreak
\definition{Definition 2.1 {\rm \cite{19, 8, 20}}}
Let $X$ be a topological space. The
{\it (Markov) free abelian topological group
over} $X$ is a pair consisting of an abelian topological group $A(X)$
and a topological embedding $X\hookrightarrow A(X)$ such that
every continuous mapping $f$ from $X$ to an abelian
topological group $G$  extends uniquely to a continuous homomorphism
$\bar f: A(X)\to G$.
\qed\enddefinition
If $X$ is a completely regular topological space
then the free abelian topological group $A(X)$  exists
and is algebraically free over the set $X$ \cite{19, 8, 20}.
A topological space $X$ is called a $k_\omega$-space \cite{18, 13-17}
if there exists a so-called $k_\omega$-decomposition
$X=\cup_{n\in\N}X_n$, where
all $X_n$ are compact, $X_n\subset X_{n+1}$ for $n\in\N$, and a subset
$A\subset X$ is closed if and only if  all intersections $A\cap
X_n,~n\in\N$, are closed.
All locally convex spaces (LCS) in this paper are real.
\definition{Definition  2.2 {\rm \cite{19, 1, 31, 6, 7, 34}}} Let $X$
be a topological space. The
{\it free locally convex space
over} $X$ is a pair consisting of a locally convex space $L(X)$ and a
topological embedding $X\hookrightarrow L(X)$ such that
every continuous mapping $f$ from $X$ to a
locally convex space $E$
extends uniquely to a continuous linear operator $\bar f: L(X)\to E$.
\qed\enddefinition
If $X$ is a completely regular topological space
then the free
locally convex space $L(X)$ exists;
the set $X$ forms a Hamel basis
for $L(X)$ \cite{31, 6, 7, 34}.
The identity mapping $id_X:X\to X$ extends to
a canonical continuous homomorphism $i:A(X)\to L(X)$.
\proclaim{Theorem 2.3 {\rm \cite{33}}}
The canonical homomorphism
$i: A(X)\hookrightarrow L(X)$ is
an embedding of $A(X)$ into the additive topological group of the LCS
$L(X)$ as a closed additive topological subgroup. \qed\endproclaim
In what follows, we will often identify $A(X)$ with a subgroup of
$L(X)$ in the above canonical way.
Denote by $L_p(X)$
the free locally convex space $L(X)$
endowed with the weak topology.
\proclaim{Theorem 2.4 {\rm \cite{6, 7}}}
Let $X$ be a completely regular space.
The canonical mapping $X\hookrightarrow L_p(X)$ is a topological
embedding, and every continuous mapping $f$ from $X$ to a
locally convex space $E$ with the weak topology
extends uniquely to a continuous linear operator $\bar f: L_p(X)\to E$.
\qed\endproclaim
The weak dual space to $L(X)$ is canonically isomorphic to
the space $C_p(X)$ of all continuous real-valued functions on $X$ with
the topology of pointwise (simple) convergence.
The spaces $L_p(X)$ and $C_p(X)$ are in duality.
Denote by $C_k(X)$ the space of continuous functions endowed with the
compact-open topology.
A topological space $X$ is called {\it Dieudonn\'e complete} \cite{5}
if its topology is induced by a complete uniformity.
For example, every Lindel\"of space  is Dieudonn\'e complete.  In
particular every  $k_\omega$-space is Dieudonn\'e complete.
\proclaim{Theorem 2.5 {\rm (Arhangel'ski\u\i\ \cite{3})}}
Let $X$ and $Y$ be Dieudonn\'e complete spaces.
If a linear mapping $C_p(X)\to C_p(Y)$ is continuous then it is
continuous as a mapping $C_k(X)\to C_k(Y)$. \qed\endproclaim
\noindent
The space $L(X)$ admits a canonical continuous monomorphism
$$L(X)\hookrightarrow C_k(C_k(X))$$
\proclaim{Theorem 2.6 {\rm (Flood \cite{6, 7}, Uspenski\u\i\
\cite{34})}}
If $X$ is a $k$-space then the monomorphism $L(X)\hookrightarrow
C_k(C_k(X))$ is an embedding of locally convex spaces.
\qed\endproclaim
Let $X$ be a topological space. A collection of continuous functions
$h_1, \dots, h_m$ on $X$ assuming their values in the closed unit
unterval $I=[0,1]$ is called
{\it basic} \cite{26, 32}
if every real-valued continuous function $f$ on $X$  can be represented
as a sum $\sum_{i=1}^n g_i\circ h_i$ of compositions of basic functions
with some continuous functions $g_i \in C (I)$.
\proclaim{2.7. Kolmogorov's Superposition Theorem {\rm \cite{11}}} The
finite-dimensional cube $I^n$ has a finite basic
family of continuous real-valued functions. \qed\endproclaim
Let us recall that for compact metrizable spaces
all three main concepts of dimension
(the covering, the small inductive and the large inductive ones)
coincide \cite{5}.
The following result is of crucial importance for us;
it is an immediate corollary
of the Kolmogorov's
Superposition Theorem,
the Menger-N\"obeling Theorem on embeddability of separable metric
spaces of dimension $\leq n$ into $\R^{2n+1}$,  and the Tietze-Urysohn
Extension Theorem \cite{5}.
\proclaim{Corollary 2.8 {\rm (Ostrand \cite{26})}}
Let $X$ be a finite-dimensional compact metrizable space. Then
there exists a finite basic family of
continuous functions on $X$.
\qed\endproclaim
For an exact upper bound on the cardinality of a
basic family of
continuous functions on a space $X$ of dimension $n$,
 see \cite{32};
however, we do not need it.

\smallpagebreak
\heading
\S 3. Auxiliary constructions
\endheading
\smallpagebreak
\proclaim{Lemma 3.1}
Consider a commutative diagram of Banach spaces and surjective
continuous linear mappings:
$$
\CD
E_1 @<r_1<< E_2 @<r_2<< E_3 @<r_3<< \dots @<r_{n-1}<<  E_n @<r_n<<
\dots\\ @V\pi_1VV @V\pi_2VV @V\pi_3VV @. @V\pi_nVV @. \\
F_1 @<q_1<< F_2 @<q_2<< F_3 @<q_3<< \dots @<q_{n-1}<< F_n @<q_n<< \dots
\endCD
$$
Denote by $E=\varprojlim E_n$ and $F=\varprojlim F_n$
the Fr\'echet spaces projective limits of corresponding inverse
sequences, and by $\pi : E\to F$ the projective limit of the mappings
$\pi_n,~n\in\N$.
Then every compact subspace $K\subset F$ is an image under the mapping
$\pi$ of a compact subspace of $E$. \endproclaim
\demo{Proof}
Let $K$ be a compact subspace of $F$. Let $K_n=q_n(K)$ for all
$n\in\N$. According to the Michael Selection Theorem (Th\. 1.4.9 in
\cite{36}),
there exists a compact subspace $C_1\subset E_1$ such that $\pi_1(C_1)=
K_1$.
Assume now that
for all $k\leq n$ we have chosen
compact subspaces $C_k\subset E_k$ such that $\pi_k(C_k)= K_k$ and
$r_{k-1}(C_k)=C_{k-1}$.
Consider the mapping
$<r_n, \pi_{n+1}>:x\mapsto (r_n(x), \pi_{n+1}(x))$
from $E_{n+1}$ to $E_n\times F_{n+1}$.
The subset
$Q_n=\{(y,z): y\in C_n, z\in K_{n+1}, q_n(z)=\pi_{n}(y)\}$
of the space $E_n\times F_{n+1}$ is compact, and is contained in the
Banach space image of the continuous linear mapping $<r_n, \pi_{n+1}>$.
Therefore, by the Michael Selection Theorem, there exists
a compact subset $C_{n+1}\subset E_{n+1}$ such that
$<r_n, \pi_{n+1}>(C_{n+1})=Q_n$.
Consequently, $r_n(C_{n+1})=C_n$, and $q_{n+1}(C_{n+1})=K_{n+1}$, which
completes the recursion step.
Finally, put $C=\varprojlim C_n$; this subset of $E$
is compact, and the property $K\subset\varprojlim K_n$ implies that
$\pi (C)=K$.
\qed\enddemo
\proclaim{Lemma 3.2}
Let $X$ and $Y$  be
$k_\omega$-spaces.
Let $h: L_p(X)\to L_p(Y)$ be an embedding of locally convex spaces.
Then $h$ is an embedding of locally convex
space $L(X)$ into $L(Y)$ as well.
\endproclaim
\demo{Proof}
As a corollary of the Hahn-Banach theorem,
the dual linear map
$h^{\ast} : C_p(Y) \to C_p(X)$ to the embedding $h$
is a continuous surjective homomorphism.
Theorem 2.5 says that $h^{\ast}$
remains continuous with respect to the compact-open
topologies on both spaces, and by virtue of the Open Mapping Theorem,
$h^{\ast}: C_k(Y)\to C_k(X)$ is open.
Since for every compact subset $C\subset X$ the elements
of the image $h(C)$ are contained in the linear span of a compact
subset of $Y$ \cite{3},
one can choose $k_\omega$-decompositions $X=\cup_{n=1}^\infty X_n$ and
$Y=\cup_{n=1}^\infty Y_n$ in such a way that for every $n\in\N$ one has
$h(sp~X_n)\subset sp~Y_n$.
It is easy to see that the restrictions mappings
$r_n: C_k(Y_{n+1})\to C_k(Y_n)$ and
$q_n: C_k(X_{n+1})\to C_k(X_n)$ are continuous surjections, and that
$C(Y)=\varprojlim C_k(Y_n)$ and $C(X)=\varprojlim C_k(X_n)$. Denote for
each $n\in\N$ by $\pi_n$ the restriction $h^{\ast}\vert_{C_k(Y_n)}$.
The conditions of Lemma 3.1
are fulfilled, and therefore every compact subset $K\subset C_k(X)$ is
an image under the mapping $h^{\ast}$ of a suitable
compact subset of $C_k(Y)$.
Therefore, the continuous linear map $h^{\ast\ast}$
dual to $h^{\ast}$
from the space $C_k(C_k(X))$ to $C_k(C_k(Y))$ is an embedding of
$C_k(C_k(X))$ into $C_k(C_k(Y))$ as a locally convex subspace. Since
the restriction of $h^{\ast\ast}$ to $L(X)$ is $h$, the desired
statement follows from Theorem 2.6.
\qed\enddemo
\proclaim{Lemma 3.3}
Let $X$ be a compact space and let $Y$ be a closed
subspace of $X$.
Denote by $\pi$ the quotient mapping from $X$ to $X/Y$.
Let $f_k,~k=1,\dots, n$ be continuous functions on $X$ such that their
restrictions to $Y$ form a basic family for $Y$, and let
$g_i,~i=1,\dots,m$ be a basic family of functions on $X/Y$.
Then the family of functions
$f_1,\dots, f_n, g_1\circ\pi, \dots, g_m\circ\pi$ is basic for $X$.
\endproclaim
\demo{Proof}
Let $f:X\to \R$ be a continuous function.
For a family of continuous functions
$h_1,\dots, h_n\in C(I)$,
the restriction $f\vert_Y$ is represented as
$\sum_{k=1}^n h_k\circ (f_k\vert_Y)=(\sum_{k=1}^n h_k\circ
f_k)\vert_Y$. Denote by $g: X\to \R$ the continuous function
$f-\sum_{k=1}^n h_k\circ f_k$; since the restriction
$g\vert_Y\equiv 0$,
the function $g$ factors through the mapping $\pi$; that is,
there exists a continuous function $h: X/Y\to I$ with
$g= h\circ\pi$.
For some collection $s_1,\dots, s_m$ of continuous functions
on $I$ one has
$h = \sum_{i=1}^m s_i\circ g_i$,
which means that $g =  \sum_{i=1}^m  s_i\circ g_i\circ\pi$.
Finally, one has
$$f = \sum_{k=1}^n h_k\circ f_k + \sum_{i=1}^m s_i\circ g_i\circ \pi,$$
as desired.
\qed\enddemo
A topological space $X$ is called {\it submetrizable} if
it admits a continuous one-to-one mapping into a metrizable space.
\proclaim{Lemma 3.4}
Let $X$ be a  submetrizable
$k_\omega$-space with $k_\omega$-decomposition
$X=\cup_{n\in\N}X_n$
such that every subspace $X_n$ is
finite-dimensional.
Then
there exists an embedding of
locally convex spaces
$\bar F: L_p(X)\hookrightarrow L_p(Y)$,
where $Y$ is the disjoint sum of countably many copies
of the closed unit interval $I$, such that
$\bar F(A(X))\subset A(Y)$.
\endproclaim
\demo{Proof}
Let $X=\cup_{n\in\N}X_n$ be a $k_\omega$ decomposition
of $X$ with $X_n\subset X_{n+1}$, for all $n\in\N$.
Since every $X_n,~n\in\N$ is a finite-dimensional metrizable
compact space, then for any $n\in\N$ so is the quotient space
$X_{n+1}/X_n$, and one can choose inductively,
using  Ostrand's Corollary 2.8 and Lemma 3.3,
a countable family of continuous functions
$f_{n,i},~n\in\N,
{}~i=1,\dots, k_n,~k_n\in\N$
from $X$ to $I$ such that for each $n\in\N$
the following are true:
\item{1.} the collection $f_{m,i},~~i=1,\dots,k_m,~
m=1,\dots, n,$ is basic for
$X_n$;
\item{2.} $f_{n+1,i}\vert_{X_n}\equiv 0$ for all $i=1,\dots,k_{n+1}$.
Denote the above family of functions $f_{n,i}$
by $\Cal F$, and let
$Y=\oplus_{f\in \Cal F}I_f$ be the disjoint sum of
countably many copies of the closed unit interval $I$.
For every $f\in \Cal F$ denote by $0_f$ the left endpoint of the closed
interval $I_f$ regarded as an element
of the free abelian group $A(Y)$.
Define a mapping, $F$, from $X$ to the free abelian group $A(Y)$
by letting
$$F(x)=\sum_{i=1,\dots, k_1}f_{1,i}(x) +\sum_{n\geq 2, ~i=1,\dots,
k_n}(f_{n,i}(x)-0_{n,i})$$
for each $x\in X$.
The mapping $F$ is properly defined, because
the first sum is finite and in the second
sum all but finitely many terms are vanishing in the free abelian group
$A(Y)$,
for every $x\in X$.
The restriction of $F$ to every $X_n$ is continuous if
being considered as a mapping to the free abelian
topological group $A(Y)$,
which fact follows from continuity of each mapping
$f_{m,i}: X_n\to I_{f_{m,i}}\subset Y,~m\leq n,~i=1,\dots, k_m$ and the
continuity of subtraction and addition in $A(Y)$. Therefore the mapping
$F:X\to A(Y)$ is continuous.
If being viewed as a continuous mapping from $X$ to the locally convex
space $L_p(Y)$, it extends to a continuous linear operator $\bar F:
L_p(X)\to L_p(Y)$.
Let
 $h:X\to \R$ be a continuous function. We will show that
there exists a continuous linear functional
$\bar h$ on the linear subspace
$\bar F(L_p(X))$
such that $\bar h\circ F\vert_X = h$.
It would mean that $\bar F(L_p(X))$ is isomorphic to
$L_p(X)$, as desired.
Construct recursively, making use of the properties 1 and 2 above, a
countable family of continuous functions
$h_{n,i},~i=1\dots, k_n,~ n\in\N$ from $I$ to $\R$ such that for every
$n\in\N$ and for all $x\in X_n$,
$$h(x) = \sum_{i=1,\dots, k_m,~ m\leq n}(h_{m,i}\circ f_{m,i})(x)$$
Let us recall that
$f_{n,i}\vert_{X_1}\equiv 0$ for all $n\geq 2$ and $i=1,\dots,k_n$. It
is easy to deduce inductively from this fact that for any $n\geq 2$
$$\sum_{i=1,\dots, k_n}h_{n,i}(0)=0$$
Define a continuous mapping $H$ from $Y$ to $I$ by letting
$H(y) =  h_{n,i}(y)$, if $y\in I_{n,i},~n\in\N$.
Extend $H$ to a continuous linear functional
$\bar H: L_p(Y)\to\R$ and denote its restriction to
$\bar F(L_p(X))$ by $\bar h$.
We claim that $\bar h\circ F\vert_{X}  = h$, or, which is the same,
that for every $n\in\N$ one has $\bar h\circ F\vert_{X_n}  = h$.
Indeed, for an arbitrary $x\in X_n$ one has:
$$(\bar h\circ F)(x)=\bar H (F(x))=
\bar H(\sum_{i=1,\dots, k_1}f_{1,i}(x) +\sum_{2\leq m \leq n,
{}~i=1,\dots, k_n}(f_{m,i}(x)-0_{m,i}))$$
$$=\sum_{i=1,\dots, k_1}H(f_{1,i}(x))+\sum_{2\leq m \leq n, ~i=1,\dots,
k_n}H(f_{m,i}(x)-0_{m,i})$$
$$=\sum_{i=1,\dots, k_m,~m\leq n}(h_{m,i}\circ f_{m,i})(x) -\sum_{2\leq
m \leq n,
{}~i=1,\dots, k_n}h_{m,i}(0)=h(x)-0=h(x).$$
\qed\enddemo
\smallpagebreak
\heading
\S 4. Main results
\endheading
\smallpagebreak
\proclaim{Theorem 4.1} For a completely regular space $X$ the following
are equivalent.
\item{(i)} The free abelian topological group $A(X)$
embeds into $A(I)$ as a topological subgroup. \item{(ii)}  The free
topological group $F(X)$ embeds into $F(I)$ as a topological subgroup.
\item{(iii)} $X$ is homeomorphic to
a closed topological subspace of $A(I)$. \item{(iv)} $X$ is
homeomorphic to
a closed topological subspace of $F(I)$. \item{(v)} $X$ is homeomorphic
to
a closed topological subspace of $\R^\infty$. \item{(vi)} $X$ is a
$k_\omega$-space such that every compact subspace of $X$ is metrizable
and finite-dimensional.
\item{(vii)} $X$ is a submetrizable $k_\omega$-space such that every
compact subspace of $X$ is
finite-dimensional.
\endproclaim
\demo{Proof}
$(i)\Rightarrow (iii)$: since the space $X$ is Lindel\"of (as a
subspace of  $A(I)$, see \cite{2}) and hence
Dieudonn\'e complete,
the group $A(X)$ is complete in its two-sided
uniformity \cite{33}
and therefore closed in $A(I)$;
but $X$ is closed in $A(X)$.
$(ii)\Leftrightarrow (iv)$: see \cite{14}.
$(iii)\Leftrightarrow (v)\Leftrightarrow (iv)$: follows from the result
of Zarichny\u\i\ \cite{37}: the free topological group $F(I)$ and the
free abelian topological group $A(I)$ are homeomorphic to open subsets
of $\R^\infty$.
$(v)\Rightarrow (vi)$:
the space $\R^\infty=\varinjlim\R^n$ is a $k_\omega$-space such that
every compact subspace of it is metrizable and finite-dimensional,
and this property is inherited by closed subsets.
$(vi)\Leftrightarrow (vii)$: see \cite{16}.
$(vii)\Rightarrow (i)$: Let
$X$ be a submetrizable
$k_\omega$ space such that
every compact subspace of $X$ is finite-dimensional. According to Lemma
3.4, there exists an embedding of
locally convex spaces $\bar F: L_p(X)\hookrightarrow L_p(Y)$, where $Y$
is the disjoint sum of
countably many copies of the closed unit interval $I$,
such that $\bar F(A(X))\subset A(Y)$.
By virtue of Lemma 3.2, $\bar F$ is also an embedding of locally convex
spaces $L(X)\hookrightarrow L(Y)$.
Its restriction to $A(X)$ is an embedding of topological groups
(Theorem 2.3).  Now apply Corollary 1.4.
\qed\enddemo

\proclaim{Theorem 4.2} For a
$k_\omega$-space $X$ the following conditions are equivalent.
\item{(i)} The free locally convex space $L(X)$ embeds into $L(I)$ as a
locally convex subspace. \item{(ii)} The free locally convex space with
the weak topology, $L_p(X)$, embeds
into $L_p(I)$ as a locally convex subspace. \item{(iii)} The space
$C_p(X)$ is a
quotient linear topological space of $C_p(I)$.
\item{(iv)} $X$ is a finite-dimensional metrizable compact space.
\endproclaim
\demo{Proof}
$(ii)\Leftrightarrow (iii)$: just dual forms of the same
statement about two locally convex spaces having their weak topology.
$(i)\Rightarrow (iv)$:
Suppose $X$ is a noncompact
$k_\omega$-space.
Since (ii) and by the same token (iii) hold,
then by virtue of Theorem 2.5 the Fr\'echet non-normable space $C_k(X)$
is an image of the Banach space $C_k(I)$ under
a surjective continuous linear mapping, which is open
by virtue of the Open Mapping Theorem -- a contradiction. Now the space
$X$, being compact, is contained
in the subspace $sp_n(Y)$
of $L(Y)$ formed by all words of the reduced
length $\leq n$ over $Y$ for some $n\in\N$ \cite{34}.
But the space $sp_n(Y)$ is a union of countably many closed subspaces
each of which is homeomorphic to a
subspace of the
$n$-th Tychonoff power of the space
$\R\times[X\oplus (-X)\oplus \{0\}]$
\cite{2}.
Therefore, $sp_n(Y)$ is finite-dimensional.
Finally, submetrizability of $X$ follows from the same property of
$L(I)$ (the latter space admits a continuous one-to-one
isomorphism into the free Banach space over $I$, \cite{1, 6, 7}).
$(iv)\Rightarrow (ii)$: it follows from Lemmas 3.4 and 3.3
that $L_p(X)$ embeds
as a locally convex subspace into the free locally convex space in the
weak topology
over a disjoint sum of finitely many homeomorphic copies of the closed
interval. The latter LCS naturally
embeds into $L_p(I)$.
$(ii)\Rightarrow (i)$: apply Lemma 3.3.
\qed\enddemo
\definition{Remark 4.3}
Surprising as it may seem, the free locally convex
space $L(\R)$ does not embed into $L(I)$, in spite of the existence of
canonical embeddings
 $A(\R)\hookrightarrow L(\R)$
and $A(I)\hookrightarrow L(I)$ and a (non-canonical one)
$A(\R)\hookrightarrow A(I)$.
It is another illustration to the well-known fact that not every
continuous homomorphism to the additive
group of reals
from a closed additive subgroup of
an (even normable) LCS extends to a continuous linear functional on the
whole space.
Such a misbehaviour is also to blame --- at least partly
--- for apparent lack of progress in attempts to make
the Pontryagin-van Kampen
duality work for free abelian topological groups \cite{24, 30}.
\enddefinition
\definition{Remark 4.4}
The problem of characterization of covering
dimension of a completely regular
space $X$ in terms of the linear topological structure of
the space $C_p(X)$ still remains open (cf. \cite{9}). However, now we
can
describe those metrizable compact spaces
having finite dimension.
\enddefinition
\proclaim{Corollary 4.5}
A metrizable compact space $X$ is finite-dimensional if and only if the
space $C_p(X)$ is a
quotient linear topological space of $C_p(I)$. \qed\endproclaim
Perhaps, the dimension of $X$ can be described in terms referring to
the linear topological structure of the space $C_p(X)$ with the help of
a characterization of dimension
in the language of basic functions due to Sternfeld \cite{32}.
\definition{Remark 4.6}
Our results also provide answers to three
problems from the
book {\it Open Problems in Topology} \cite{27}.
\enddefinition
\item{}{\smc Problem 511.} {\sl
Is $A(I^2)$ topologically isomorphic with a subgroup of $A(I)$?}
\smallskip
\itemitem{}{\bf Yes} (cf. Theorem 4.1).
\smallskip
\item{}{\smc Problem 1046.} {\sl
 Assume that $C_p(X)$ can be mapped by a linear continuous
           mapping onto $C_p(Y)$. Is it true that $\dim Y\leq\dim X$ ?
           What if $X$ and $Y$ are compact?}

\item{}{\smc Problem 1047.} {\sl
 Assume that $C_p(X)$ can be mapped by an open linear
           continuous mapping onto $C_p(Y)$. Is it true that
           $\dim Y\leq\dim X$ ? What if $X$ and $Y$ are compact?}
\smallskip
\itemitem{}{\bf No}, in all four cases (cf. Corollary 4.5).

\smallpagebreak
\heading Acknowledgment
\endheading
\smallpagebreak
The third author (V.G.P.) thanks
the University of Wollongong for hospitality extended in
August-September 1992.

\bye